\documentclass[9pt]{amsart}
\usepackage{amsmath,amsthm,amsfonts,amssymb,amscd,euscript}
\usepackage{color}
\usepackage{graphics}
\usepackage{showlabels}
\input{xy}
\xyoption{all}
\newlength{\defbaselineskip}
\theoremstyle{plain}

\theoremstyle{definition}

\begin{document}

\title{Every finite acyclic quiver is a full subquiver of a quiver mutation equivalent to a bipartite quiver}
\author{Kyungyong Lee}
\address{Department of Mathematics, Wayne State University, Detroit, MI 48202}
\email{klee@math.wayne.edu}
\thanks{{The author has been  partially  supported by the NSF grant DMS 0901367. }}

\subjclass[2010]{13F60}
\date{\today} 

\keywords{cluster algebra, positivity conjecture}
\maketitle
We give a very short proof of the claim in the title. Under certain circumstances, this result reduces the study of acyclic cluster algebras to that of cluster algebras associated to bipartite quivers. For instance, it implies that the positivity conjecture for skew-symmetric acyclic cluster algebras, whose known proofs \cite{KQ,LS4,DMSS} are quite technical,  is an immediate consequence of a result of Nakajima in \cite{N}. 

A quiver $Q = (Q_0, Q_1, t, h)$ is a finite oriented graph consisting of $Q_0$, the set of vertices, $Q_1$, the set of arrows, and two functions $t, h : Q_1 \longrightarrow Q_0$ attaching to an arrow $\alpha$ its tail $t(\alpha)$ and head $h(\alpha)$. The quiver mutation has been discovered in \cite{FZ}. Let $Q$ be any quiver without oriented cycles. Let $\ell$ be the maximum length of all oriented paths in $Q$. Let $Q^{(0)}=Q$.   For each integer $i\geq 0$, we construct a new quiver $Q^{(i+1)}$ from $Q^{(i)}$ as follows. 

\noindent\textbf{Step 1:}  Let $\Gamma_i$ be the set of oriented paths  in $Q^{(i)}$ with length $\ell$. Let $\alpha\in Q_1^{(i)}$ be an arrow which does not lie in any  oriented path with length $\ell$, i.e., $\alpha\not\in \gamma$ for all $\gamma\in\Gamma_i$. (If no such arrow exists, then let $j=i$ and  go to Step 2.)  Introduce a new vertex $v_{\alpha}$. Add an arrow from $h(\alpha)$ to $v_{\alpha}$ and another arrow from $v_{\alpha}$ to $t(\alpha)$. Let $Q^{(i+1)}$ be the resulting quiver obtained by mutating at $v_{\alpha}$. Note that $Q^{(i+1)}$ does not have oriented cycles. Repeat Step 1. This step terminates in finitely many steps.

\noindent\textbf{Step 2:} Note that all maximal oriented paths in $Q^{(j)}$ have the same length. For each integer $i\geq j$, let $Q^{(i+1)}$ be the quiver obtained from  $Q^{(i)}$ by mutating  at all sources in  $Q^{(i)}$. Then the length of every maximal oriented path in $Q^{(i)}$ is equal to 1 or $\max(1,\ell+j-i)$ for $i\ge j$. In particular, $Q^{(i)}$ is bipartite for $i\geq j+\ell-1$. \\
\noindent\textbf{Example.}
{\tiny
$
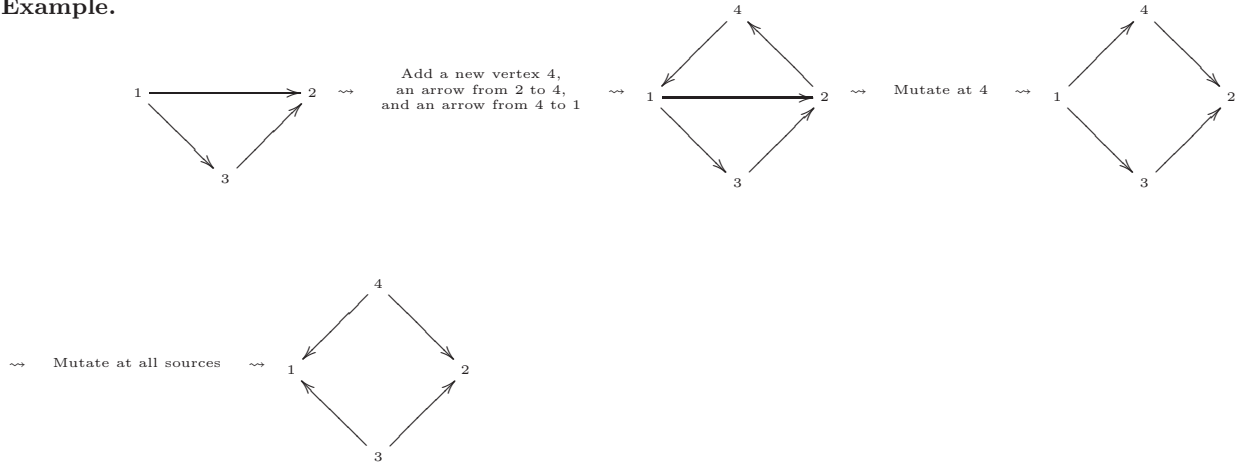
\xymatrix{& \\ 1\ar[rr] \ar[rd] &&2\\&3\ar[ru]} $    	$\xymatrix{\\	 \leadsto \\ }$ $\begin{array}{c}\\ \\ \\ \\ \\ \\ \\ \\ \\ \\ \text{Add a new vertex }4,\\ 
\text{an arrow from 2 to 4,}\\ \text{and an arrow from 4 to 1}\end{array}$   $\xymatrix{\\	 \leadsto \\ }$   $\xymatrix{  &4 \ar[ld] & \\ 1\ar[rr] \ar[rd] & &2\ar[lu]\\  &3\ar[ru]}$ $\xymatrix{\\	 \leadsto \\ }$ $\begin{array}{c}\\ \\ \\ \\ \\ \\ \\ \\ \\ \\ \text{Mutate at }4\\ 
\end{array}$   $\xymatrix{\\	 \leadsto \\ }$
$\xymatrix{  &4 \ar[rd] & \\ 1\ar[ru] \ar[rd] & &2\\  &3\ar[ru]}$ $\xymatrix{\\	 \leadsto \\ }$ $\begin{array}{c}\\ \\ \\ \\ \\ \\ \\ \\ \\ \\ \text{Mutate at all sources}\\ 
\end{array}$   $\xymatrix{\\	 \leadsto \\ }$
$\xymatrix{  &4 \ar[rd]\ar[ld]  & \\ 1 & &2\\  &3\ar[ru]\ar[lu]}$ 
}

\noindent \textbf{Remark.} Applying the idea of embedding (followed  by mutations) to general quivers, one can show that  every quiver without 2-cycles is a full subquiver of a quiver mutation equivalent to another quiver whose simple cycles have the same positive weight. 

\noindent\emph{Acknowledgements.} We are grateful to Kyu-Hwan Lee,   Hiraku Nakajima and Fan Qin for valuable discussions.


\begin{thebibliography}{99}
\bibitem{DMSS}
B.~Davison, D.~Maulik, J.~Sch\"urmann and B.~Szendr\H{o}i, Purity for graded potentials and quantum cluster positivity, arXiv:1307.3379.

\bibitem{FZ}
S.~Fomin and  A.~Zelevinsky, Cluster algebras I: Foundations, J. Amer. Math. Soc.  \textbf{15} (2002) No. 2  497--529.  


\bibitem{KQ}
Y.~Kimura and F.~Qin, Graded quiver varieties, quantum cluster algebras, and dual canonical basis, arXiv:1205.2066v2.

\bibitem{LS4}
K.~Lee and R.~Schiffler, Positivity for cluster algebras, arXiv:1306.2415.

\bibitem{N}
H.~Nakajima, Quiver varieties and cluster algebras, Kyoto J. Math. \textbf{51} (2011), No. 1, 71--126. 
\end{thebibliography}
\end{document}